\def \a{\alpha}
\def \b{\beta}
\def \g{\gamma}
\def \d{\delta}

\def \lam{\lambda}

\def \r{\right}

\def \be{\begin{equation}}
\def \eq{\end{equation}}
\def \bee{\begin{eqnarray}}
\def \eqq{\end{eqnarray}}

\def \zb{\bar{z}}
\def \wb{\bar{w}}
\def \del{\partial}
\def \delb{\bar{\partial}}
\def \Zp {{\cal Z_+}}
\def \Zm {{\cal Z_-}}
\def \H {{\cal H}}

\def \O  {{\cal O}}
\def \wb {\bar{w}}

\def \S {$S_q^2$ }
\def \SU{$SU_q (2)$ }
\def \U {$\cal U $ }

\def \deltb {\bar{\delta}}

\def \r {\rho}

\documentstyle[12pt]{article}

\begin{document}
\begin{titlepage}
\begin{center}
March 29, 1995
        \hfill  LBL-37055\\
          \hfill    UCB-PTH-95/10 \\
\vskip .5in

{\LARGE } The quantum 2-sphere as a complex quantum manifold\footnote{This work
was supported in part by the Director, Office of
Energy Research, Office of High Energy and Nuclear Physics, Division of
High Energy Physics of the U.S. Department of Energy under Contract
DE-AC03-76SF00098 and in part by the National Science Foundation under
grant PHY-90-21139.}
\vskip .5in
{Chong-Sun Chu, Pei-Ming Ho and Bruno Zumino}

\vskip .2in
{\em
Department of Physics \\University of California \\ and\\
   Theoretical Physics Group\\
    Lawrence Berkeley Laboratory\\
      University of California\\
    Berkeley, CA 94720}
\end{center}

\vskip .5in

\begin{abstract}

We describe the quantum sphere of Podle\'{s} for $c=0$
by means of a stereographic projection which is analogous to that
which exibits the classical sphere as a complex manifold.
We show that the algebra of functions and the differential calculus
on the sphere are covariant under the coaction of fractional transformations
with \SU coefficients as well as under the action of \SU vector fields.
Going to the classical limit we obtain the Poisson sphere.
Finally, we study the invariant integration of functions on the sphere
and find its relation with the translationally invariant integration
on the complex quantum plane.

\end{abstract}

\end{titlepage}
\renewcommand{\thepage}{\roman{page}}
\setcounter{page}{2}
\mbox{ }

\vskip 1in

\begin{center}
{\bf Disclaimer}
\end{center}
\vskip .2in
\begin{scriptsize}
\begin{quotation}
This document was prepared as an account of work sponsored by the United
States Government. While this document is believed to contain correct
 information, neither the United States Government nor any agency
thereof, nor The Regents of the University of California, nor any of their
employees, makes any warranty, express or implied, or assumes any legal
liability or responsibility for the accuracy, completeness, or usefulness
of any information, apparatus, product, or process disclosed, or represents
that its use would not infringe privately owned rights.  Reference herein
to any specific commercial products process, or service by its trade name,
trademark, manufacturer, or otherwise, does not necessarily constitute or
imply its endorsement, recommendation, or favoring by the United States
Government or any agency thereof, or The Regents of the University of
California.  The views and opinions of authors expressed herein do not
necessarily state or reflect those of the United States Government or any
agency thereof of The Regents of the University of California and shall
not be used for advertising or product endorsement purposes.
\end{quotation}
\end{scriptsize}

\vskip 2in

\begin{center}
\begin{small}
{\it Lawrence Berkeley Laboratory is an equal opportunity employer.}
\end{small}
\end{center}
\newpage
\renewcommand{\thepage}{\arabic{page}}
\setcounter{page}{1}


\section{INTRODUCTION}

Quantum spheres can be defined in any number of dimensions by
normalizing a vector of quantum Euclidean space\cite{FRT}.
The differential calculus on quantum Euclidean space\cite{OZ}
induces a calculus on the quantum sphere.
The case of two-spheres in three space is special in that
there are many more possibilities than the one obtained
from the general construction.
These have been studied by P. Podle\'{s}\cite{P1,P2,P3,P4}
who has also shown how to define a noncommutative differential calculus on
them.
In this paper we study in detail a particular case of Podle\'{s} spheres
which is one of those special to three space dimensions.
In this case the algebra of functions on the sphere is a subalgebra
of the algebra of functions on \SU and the differential calculus
on the sphere can be inferred from a differential calculus on \SU.
We can also define a stereographic projection and describe the coaction of \SU
on the sphere by fractional transformations on the complex variable in the
plane
analogous to the classical ones.
The quantum sphere appears then as the quantum deformation
of the classical two-sphere described as a complex manifold.

Our quantization of the sphere is not symmetric between the north and the south
pole.
This asymmetry is also apparent when we go to the classical limit of the
Poisson sphere
\cite{We} and it seems to be unavoidable in our approach.
A description of Podle\'{s} spheres was given in an interesting paper by
\v{S}\v{t}ov\'{i}\v{c}ek\cite{St}.
He shows that the sphere can be understood as the patching of two complex
quantum planes.
His choice of variable is symmetric between the two planes,
but the coaction of \SU is very complicated in terms of his variable.
Also, \v{S}\v{t}ov\'{i}\v{c}ek does not consider the noncommutative calculus on
the sphere.

\section{\S AS A COMPLEX MANIFOLD}

In Ref.\cite{P1}, a family of quantum 2-spheres was introduced. There, the
algebra of functions over the sphere is generated by 3 coordinates, subjected
to
a condition that reduces the number of indepedent generators to 2.  The case of
$c=0$ is of special interest\cite{BM}. In this case, the algebra is generated
by $b_+=\g \d, b_-=\a \b,  b_3=\a \d$, (where $\a, \b,\g,\d \in SU_q(2)$) with
commutations
\bee
   &b_3b_-=(1-q^{-2})b_-+q^{-2} b_-b_3,\\
   &b_3b_+=b_+(1-q^{2})+q^{2}b_+b_3, \\
   &q^{-2}b_-b_+=q^{2}b_+b_- + (q^{-1} -
    q)(b_3-1),
\eqq
and constraint
\be b_3^2=b_3+q^{-1}b_-b_+. \eq
The $*$-algebra structure is $b_\pm^*=-q^{\mp 1} b_\mp$,
$b_3^*=b_3$, and $q^{*} = q$.

One can construct a stereographic projection to go from the 3 coordinates
$b_\pm, b_3$ to the complex plane $z, \zb$. Define
\bee \label{stereo}
   &z=-q b_- (1-b_3)^{-1}=\a \g^{-1}, \label{SU1}\\
   &\zb = b_+ (1-b_3)^{-1}= -\d \b^{-1}, \label{SU2}
\eqq
which is the projection from the north pole of the sphere to the plane with
coordinates $z, \zb$.
It is easy to derive the commutation relation
\be \label {zz}
 z \zb =q^{-2} \zb z +q^{-2}-1 \eq
and the $*$-structure $z^*=\zb$.
This differs from the usual quantum plane by an additional inhomogeneous
constant term. One can check directly that Eq.(\ref {zz}) is covariant under
the fractional transformation,
with $\pmatrix{a&b\cr c&d} \in SU_q(2)$,
\be \label{z-transf}
z \rightarrow (a z +b)(c z+d)^{-1}, \quad
\zb \rightarrow -(c-d \zb)(a-b \zb)^{-1},
\eq
which is induced from the \SU coproduct, interpreted as a left transformation.
Here $a,b,c$ and $d$ commute with $z$ and $\zb$.

\section{DIFFERENTIAL CALCULUS}

In Refs.\cite{P2,P3,P4}, differential structures on \S are studied and
classified.
In this section we give a differential calculus on \S in terms of the complex
coordinates $z$ and $\zb$.
Just as the algebras of functions and vector fields on \S can be inferred from
those of \SU,
so can the differential calculus.

For \SU there are several well-known calculi\cite{Wcal1,Wcal2}: the 3D left-
and right-covariant differential calculi, and the $4D_{+}, 4D_{-}$ bi-covariant
calculi.
The 4D bi-covariant calculi have one extra dimension in their space of
one-forms compared with the classical case.
The right-covariant calculus will not give a projection on \S in a closed form
in terms of $z$, $\zb$, which are defined to transform from the left.
Therefore we shall choose the left-covariant differential calculus.

It is straightforward to obtain the following relations from those for \SU:
\bee
   &z dz=q^{-2}dz z,     & \zb dz=q^{2}dz \zb,  \label{zdz}\\
   &z d\zb=q^{-2}d\zb z, & \zb d\zb=q^{2}d\zb \zb, \\
   &(dz)^{2}=(d\zb)^{2}=0, \label{dzdz}
\eqq
and
\bee
   &dz d\zb=-q^{-2}d\zb dz.
\eqq

We can also define derivatives $\del$, $\delb$ such that on functions,
\be
d=dz\del+d\zb\delb.
\eq
{}From the requirement $d^{2}=0$ and the undeformed Leibniz rule for $d$
together with
Eqs. (\ref{zdz}) to (\ref{dzdz}) it follows that:
\begin{eqnarray}
\label{delz}
   &\del z=1+q^{-2}z \del, &\del \zb=q^{2}\zb \del,   \\
   &\delb z=q^{-2}z \delb, &\delb \zb=1+q^{2}\zb \delb,
\eqq
and
\bee
\label{deldelb}
   &\del \delb=q^{-2}\delb \del.
\eqq
It can be checked explicitely that these commutation relations are covariant
under the transformation (\ref{z-transf}) and
\bee
   &dz\rightarrow dz(q^{-1}cz+d)^{-1}(cz+d)^{-1}, \\
   &\del\rightarrow (cz+d)(q^{-1}cz+d)\del,
\eqq
which follow from (\ref{z-transf}) and the fact that $d$ is invariant.

The $*$-structure also follows from that of \SU:
\begin{eqnarray}
   &(dz)^{*}=d\zb,  \label{dz*}\\
   &\del^{*}=-q^{-2}\delb+(1+q^{-2})z \rho^{-1},  \label{del*}\\
   &\delb^{*}=-q^{2}\del+(1+q^{2}) \rho^{-1}\zb. \label{delb*}
\end{eqnarray}
where we have introduced
\be \label{rho}
   \rho=1+\zb z
\eq
(remember that the $*$-involution inverts the order of factors in a product).

The inhomogeneous pieces on the RHS of the Eqs.(\ref{del*}) and (\ref{delb*})
reflect the fact that
the sphere has curvature.
Incidentally all the commutation relations in this section admit another
possible involution:
\bee
\label{Cdz*}     &(dz)^{*}=d\zb, \\
\label {Cdel*}   &\del^{*}=-q^{2}\delb, \\
\label{Cdelb*}   &\delb^{*}=-q^{-2}\del.
\eqq
This involution is not covariant under the fractional transformations and
cannot be used for the sphere.
However, it can be used when we have a quantum plane defined by the same
algebra of functions and calculus.

We shall take Eqs. (\ref{zdz}) to (\ref{rho}) as the definition of the
differential calculus on \S.

It is interesting to note that there exist two different types of symmetries in
the calculus.
The first symmetry is that if we put a bar on all unbarred variables ($z$,
$dz$, $\del$),
take away the bar from any barred ones and at the same time replace $q$ by
$1/q$ in any statement about the calculus, the statement is still true.

The second symmetry is the consecutive operation of the two $*$-involutions
above, so that
\bee
   &\del\rightarrow -q^{2}\delb^{*}=q^{4}\del-q^{2}(1+q^{2}) \rho^{-1}\zb, \\
   &\delb\rightarrow -q^{-2}\del^{*}=q^{-4}\delb-q^{-2}(1+q^{-2})z \rho^{-1},
\eqq
with $z, \zb, dz, d\zb$  unchanged.
This replacement can be iterated $n$ times and gives a symmetry which resembles
that of a gauge transformation on a line bundle:
\bee
   \del\rightarrow \del^{(n)} &\equiv& q^{4n}\del-q^{2}[2n]_{q} \rho^{-1}\zb \\
                              &=& q^{4n} \rho^{2n}\del \rho^{-2n}, \\
   \delb\rightarrow \delb^{(n)} &\equiv& q^{-4n}\delb-q^{-2}[2n]_{1/q}z
\rho^{-1} \\
                                &=& q^{-4n} \rho^{2n}\delb \rho^{-2n},
\eqq
where $[n]_{q}=\frac{q^{2n}-1}{q^{2}-1}$.
For example, we have
\be \del^{(n)} z = 1 + q^{-2} z \del^{(n)}.
\eq
Making a particular choice of $\del, \delb$ is like fixing a gauge.

Many of the features of a calculus on a classical complex manifold are
preserved.
Define $\delta=dz\del$ and $\deltb=d\zb\delb$ as the exterior derivatives on
the holomorphic and antiholomorphic functions on \S respectively. We have:
\bee
   &\left[\delta,z\right]=dz, \quad \left[\delta,\zb\right]=0, \\
   &\left[\deltb,z\right]=0, \quad \left[\deltb,\zb\right]=d\zb, \\
   &d=\delta +\deltb.
\eqq
The action of $\delta$ and $\deltb$ can be extended consistently on forms as
follows
\bee
   &\delta dz=dz \delta=0, \quad \deltb d\zb=d\zb \deltb=0,\\
   &\{\delta, d\zb\}=0, \quad \{\deltb, dz\}=0,\\
   &\delta^{2}=\deltb^{2}=0,\\
   &\{\delta, \deltb\}=0,
\eqq
where $\{\cdot,\cdot\}$, $[\cdot,\cdot]$ are the anticommutator and commutator
respectively.

\section{THE RIGHT INVARIANT VECTOR FIELDS ON \S}

In this section we want to define vector fields on \S which generate the
fractional transformation mentioned above.
We will see that these vector fields can be inferred from those on \SU.

First let us recall some well-known facts about the vector fields on \SU (see
for example Ref.\cite{Z}).
The enveloping algebra \U of \SU is usually said to be generated by the
left-invariant vector fields $H_{L}$,$X_{L\pm}$
which are arranged in two matrices $L^+$ and $L^-$.
The action of these vector fields corresponds to infinitesimal right
transformation: $T\rightarrow TT'$.
What we want now is the infinitesimal version of the left transformation given
by Eqs.(\ref {z-transf}),
hence we shall use the right-invariant vector fields $H_{R}$,$X_{R\pm}$.
Since only the right-invariant ones will be used, we will drop the subscript
${R}$ hereafter.

The properties of the right-invariant vector fields are similar to those of the
left-invariant ones.
Note that if an \SU matrix $T$ is transformed from the right by another \SU
matrix $T'$,
then it is equivalent to say that the $SU_{1/q}(2)$ matrix $T^{-1}$ is
transformed
from the left by another $SU_{1/q}(2)$ matrix $T'^{-1}$.
Therefore one can simply write down all properties of the left-invariant vector
fields and then make the replacements: $q\rightarrow 1/q$, $T\rightarrow
T^{-1}$
and left-invariant fields$\rightarrow$right-invariant fields.

Using the matrices:
\bee
   M^{+}= \left(
            \begin{array}{ll}
              q^{-H/2} & q^{-1/2}\lam X_{+} \\
              0        & q^{H/2}
            \end{array}
          \right),
 &
   M^{-}= \left(
            \begin{array}{ll}
              q^{H/2}         & 0        \\
              -q^{1/2}\lam X_{-} & q^{-H/2}
            \end{array}
          \right),
\eqq
the commutation relations between the vector fields are given by,
\bee
   R_{12}M_{2}^{+}M_{1}^{+}=M_{1}^{+}M_{2}^{+}R_{12}, \\
   R_{12}M_{2}^{-}M_{1}^{-}=M_{1}^{-}M_{2}^{-}R_{12}, \\
   R_{12}M_{2}^{+}M_{1}^{-}=M_{1}^{-}M_{2}^{+}R_{12},
\eqq
while the commutation relations between the vector fields and the elements of
the quantum
matrix in the smash product of \U and \SU are,
\bee
   T_{1}M_{2}^{+}=M_{2}^{+}{\cal R}_{12}T_{1}, \\
   T_{1}M_{2}^{-}=M_{2}^{-}{\cal R}_{21}^{-1}T_{1},
\eqq
where $T$ is a \SU matrix, ${\cal R}=q^{-1/2}R$ and $R$ is the $GL_q(2)$
R-matrix.
Clearly $M^{+}$, and $M^{-}$ are the right-invariant counterparts of $L^{+}$
and $L^{-}$.
The commutation relations between the $M$'s and the $T$'s tell us how the
functions on \SU are
transformed by the vector fields $H$,$X_{+}$,$X_{-}$.
It is convenient to define a different basis for the vector fields,
\bee
  &\Zp=X_{+} q^{H/2}, \\
  &\Zm=q^{H/2}X_{-}
\eqq
and
\bee
   &\H=[H]_{q}=\frac{q^{2H}-1}{q^{2}-1}.
\eqq
They satisfy the commutation relations
\bee \label{ZZH}
&\H \Zp - q^4 \Zp \H =(1+q^2)\Zp, \\
    &\Zm \H - q^4 \H \Zm =(1+q^2) \Zm
\eqq
and
\be
   q\Zp\Zm-q^{-1}\Zm\Zp=\H.
\eq
Using the expressions of $z, \zb$ in terms of $\a, \b, \g, \d$, one can easily
find the action
of the vector fields on the variables $z, \zb$ on the sphere,
\bee
   &\Zp z=q^{2}z\Zp+q^{1/2}z^{2}, \label{vect1}\\
   &\Zp\zb=q^{-2}\zb\Zp +q^{-3/2}, \label{Zpzb} \\
   &\H z=q^{4}z\H+(1+q^{2})z, \label{Hz} \\
   &\H\zb=q^{-4}\zb\H-q^{-4}(1+q^{2})\zb, \label{Hzb} \\
   &\Zm z =q^2 z \Zm -q^{1/2} \label{Zmz}
\eqq
and
\bee
   &\Zm \zb= q^{-2} \zb \Zm -q^{-3/2} \zb^2. \label{vect2}
\eqq
It is clear that a $*$-involution can be given:
\be \label{Z*}
\Zp^*=\Zm, \quad  \H^*=\H.
\eq

Since all the relations listed above are closed in the vector fields and $z$,
$\zb$
(this would not be the case if we had used the left-invariant fields),
we can now take these equations as the definition of the vector fields that
generate the fractional tranformation on \S.
We shall take our vector fields to commute with the exterior differentiation
$d$.
This is consistent for right-invariant vector fields in a left-covariant
calculus
and allows us to obtain the action of our vector fields on the differentials
$dz$ and $d\zb$,
as well as on the derivatives $\del$ and $\delb$.
For instance (\ref{vect1}) gives
\be
   \Zp dz = q^{2}dz\Zp + q^{1/2}(dzz+zdz)
\eq
and
\be
   \del\Zp = q^{2}\Zp\del + q^{-3/2}(1+q^{2})z\del.
\eq

\section{ MORE ABOUT THE CALCULUS }

The calculus described in the previous section has a very interesting property.
There exists a one-form $\Xi$ having the property that
\be \label{dxi}
   \Xi f\mp f\Xi = \lambda df,
\eq
where, as usual, the minus sign applies for functions or even forms
and the plus sign for odd forms.
Indeed, it is very easy to check that
\be
   \Xi = \xi - \xi^*
\eq
\be
   \xi = qdz \rho^{-1}\zb
\eq
satisfies Eq.(\ref{dxi}) and
\be
   \Xi^* = -\Xi.
\eq
It is also easy to check that
\be
   d\Xi = 2qd\zb \rho^{-2}dz
\eq
and
\be
   \Xi^{2} = q\lambda d\zb \rho^{-2}dz.
\eq
Suitably normalized, $d\Xi$ is the natural area element on the quantum sphere.
Notice that $\Xi^{2}$ commutes with all functions and forms,
as required for consistency with the relation
\be
   d^{2} = 0.
\eq

The existence of the form $\Xi$ within the algebra of
$z, \zb, dz, d\zb$ is especially interesting because no such form exists
for the 3-D calculus on $SU_{q}(2)$,
from which we have derived the calculus on the quantum sphere
(a one-form analogous to $\Xi$ does exist for the two
bicovariant calculi on $SU_{q}(2)$,
but we have explained before why we didn't choose either of them).
It is also interesting that $d\Xi$ and $\Xi^{2}$ do not vanish
(as the corresponding expressions do in the bicovariant calculi on the quantum
groups
or in the calculus on quantum Euclidean space).
We see here an example of Connes' calculus\cite{Con} of the type $F^{2} = 1$
rather than $F^{2} = 0$.

The one-form $\Xi$ is regular everywhere on the sphere,
except at the point $z = \zb =\infty$,
which classically corresponds to the north pole.
We shall discuss this question in Sec.\ref{Poisson}
where we argue that the pole singularity at that point
can be included by allowing forms with distribution valued coefficients.
The area element $d\Xi$ is regular everywhere on the sphere.

It is interesting to see how $\Xi$ and $d\Xi$ transform
under the action of the right invariant vector fields
or under the coaction of the fractional transformations (\ref{z-transf}).
Using (\ref{vect1}) to (\ref{vect2}) one finds
\be
   \Zp\Xi = \Xi\Zp + q^{-1/2}dz \label{ZpXi}
\eq
and
\be
   \H\Xi = \Xi\H. \label{HXi}
\eq
These equations are consistent with (\ref{dxi}).
For instance,
\bee
   \Zp(\lambda dz-\Xi z+z\Xi)& =& q^{2}(\lambda dz-\Xi z+z\Xi) \Zp+ \nonumber
\\
    & &  + q^{1/2}(\lambda dz^{2}-\Xi z^{2}+z^{2}\Xi)- \nonumber \\
    & &- q^{-1/2}(dz z-q^{2}zdz).
\eqq
Eqs. (\ref{ZpXi}) and (\ref{HXi}) imply that $d\Xi$ commutes with $Z_{\pm}$ and
$\H$,
as expected for the invariant area element.

For the fractional transformation (\ref{z-transf}) one finds $\xi\rightarrow
\xi'$ where
\be
   \xi'-\xi = -q(dz)cd^{-1}(1+cd^{-1}z)^{-1} \label{xi-transf}
\eq
and a similar formula for $\xi^*$.
The right hand side of (\ref{xi-transf}) is a closed one-form,
since $(dz)^{2} = 0$, so one could write
\be
   \xi'-\xi = -qd[\log_{q}(1+cd^{-1}z)]
\eq
with a suitably defined quantum function $\log_{q}$.
At any rate
\be
   d\xi' = d\xi
\eq
so that the area element two-form is invariant under
finite transformations as well.

\section{ PATCHING TWO QUANTUM PLANES }

The variables $z$ and $\zb$ cover the sphere with the exception of the north
pole. In
analogy with the classical case, we can introduce new variables $w=z^{-1}$ and
$\wb=\zb^{-1}$
which describe the sphere without the south pole. These variables satisfy the
commutation relation
\be \label {wwb}
 w \wb =q^{-2} \wb w +(q^{-2}-1)w \wb^2 w
\eq
which is covariant under the transformation
\be w \rightarrow (d w +c) (b w +a)^{-1}, \quad \wb \rightarrow -(a \wb -b)(c
\wb-d)^{-1}.
\eq
Notice that the commutation relation (\ref{wwb}) is different from that
satisfied by
$z$ and $\zb$; our way of quantizing the sphere is inherently asymmetric
between the
north and the south pole.

The calculus in $z$ and $\zb$ induces a calculus in $w$ and $\wb$. It is not
hard to
derive the commutation relations  for this $w, \wb$ calculus as well as the
mixed
commutation relations.
For example, we have
\bee
   wdw = q^{2}dww, \\
   \del_{w}w = 1+q^{2}w\del_{w}
\eqq
and
\be
   dzw = q^{-2}wdz.
\eq

Since $w$ and $\wb$ are functions of $z$ and $\zb$, Eq.(\ref{dxi}) is valid for
functions
and forms in $w$ and $\wb$, with the same $\Xi$. In terms of $w$ and $\wb$ the
one-forms
$\xi$ and $\xi^*$ are given by
\be \xi = -w^{-1} dw (1+\wb w)^{-1}, \quad \xi^* =-(1+\wb w)^{-1} d\wb
\wb^{-1}. \eq
Clearly they are singular at the north pole $w =\wb =0$. This polar singularity
is an
intrinsic feature of our asymmetric quantization and of our calculus.  We
believe that
it can be controlled by allowing distributions, rather than just functions as
the
elements of our algebra and as coefficients of differential forms. In order to
avoid
the need to develop the concept of distribution in the framework of
noncommutative
algebra, we explain our point of view in the next section for the limit of the
Poisson sphere.

\section{THE POISSON SPHERE}\label{Poisson}
The commutation relations of the previous sections give us, in the limit $q
\rightarrow 1$,
a Poisson structure on the sphere. The Poisson Brackets (P.B.s) are obtained as
usual as  a limit
\be (f,g) = \lim_{h \rightarrow 0} \frac{fg-gh}{h}, \quad q^2=e^h=1+h+[h^2].
\eq
For instance, the commutation relation (\ref{zz}) gives
\be z \zb =(1-h) \zb z -h +[h^2] \eq
and therefore
\be (\zb, z)=\rho. \eq
Similarly one finds
\bee
&(dz,z)=zdz, &(d\zb,z)=z d\zb, \\
&(dz, \zb)=-\zb dz, &(d\zb, \zb)=-\zb d\zb
\eqq
and
\bee &(d \zb, dz)=d\zb dz. \eqq
In this classical limit functions and forms commute or anticommute according to
their even
or odd parity, as usual. The P.B. of any quantity with itself vanishes. The
P.B. of two even
quantities or of an even and an odd quantity is antisymmetric, that of two odd
quantities is
symmetric. It is
\be d(f,g)=(df, g) \pm (f, dg) \eq
where the plus (minus) sign applies for even (odd) $f$. Notice that we have
enlarged the concept
of Poisson bracket to include differential forms. This is very natural when
considering the classical
limit of our commutation relations.

In the classical limit, Eq.(\ref{dxi}) becomes
\be \label{pb-dxi} (\Xi, f)=df \eq
where
\be \Xi=\xi-\xi^* \eq
and
\be
     \xi=dz \zb \rho^{-1}, \quad \xi^*=d\zb z \rho^{-1}
\eq
are ordinary classical differential forms. Now
\be d\Xi=2d\zb dz \rho^{-2} \eq
and
\be \Xi^2=0. \eq

As before, the variables $z$ and $\zb$ cover the sphere except for the north
pole, while $w$ and $\wb$ miss
the south pole. It is
\be (\wb, w)=\wb w(1+\wb w). \eq
The Poisson structure is not symmetric between the north and south pole. All
P.B.s of regular functions and forms
vanish at the north pole $w=\wb=0$. Therefore, for Eq.(\ref{pb-dxi}) to be
valid, the one-form $\Xi$ must be singular
at the north pole. Indeed one finds
\be \label{xi} \xi =\frac{dw \wb}{1+\wb w} -\frac{dw}{w}, \quad
\xi^* =\frac{d\wb w}{1+\wb w} -\frac{d\wb}{\wb},
\eq
and
\be \Xi=\frac{w d\wb-\wb dw}{\wb w (1+\wb w)}. \eq
On the other hand the area two-form
\be \label{dXi0} d \Xi =2 \frac{d\wb dw}{(1+\wb w)^2} \equiv \Omega \eq
is regular everywhere on the sphere.

The singularity of $\Xi$ at the north pole is not a real problem if we treat it
in the sense  of
the theory of distributions. Consider a circle $C$ of radius $r$ encircling the
origin of the $w$ plane in
a counter-clockwise direction and set
\bee &w=r e^{i\theta}, &\wb =r e^{-i \theta}. \eqq
Using (\ref{xi}), we have
\be
\int_C \Xi = \int_C \frac {\wb dw-w d\wb}{1+\wb w}  -4 \pi i.
\eq
As $r \rightarrow 0$ the integral in the right hand side tends to zero
because the integrand is regular at the origin.
Stokes theorem can be satisfied
even at the origin if we modify Eq.(\ref{dXi0}) to read
\be
   d\Xi= \Omega -4 \pi i \d(w) \d(\wb) d\wb dw. \label{dxi1}
\eq
It is
\be \int_{S^2} \Omega =4 \pi i \eq
so that
\be \int_{S^2} d \Xi =0 \eq
as it should be for a compact manifold without boundary.
Notice that the additional delta function term in (\ref{dxi1})
also has zero P.B.s with all functions and forms as required by consistency.

\section{INTEGRATION}

We now return to the quantum case.
For the integral of a function $f$ over the sphere
we shall use the notation $<f>$.
A left-invariant integral can be defined,
up to a normalization constant,
by requiring invariance under the action of the right-invariant
vector fields
\be
   <\O f(z,\zb)> = 0, \;\; \O = \Zp,\Zm,\H.
\eq

Using $\H$ and Eqs.(\ref{Hz}) and (\ref{Hzb}) one finds that
\be
   <z^{k}\zb^{l}g(\zb z)> = 0,\;\; unless \;\; k=l.
\eq
(Here $g$ is a convergence function.)
Therefore we can restrict ourselves to integrals of the form $<f(\zb z)>$.

Eqs.(\ref{vect1}) and (\ref{Zpzb}) imply
\be
   \Zp\r = \r\Zp + q^{1/2}z\r
\eq
and
\be
   \Zp\r^{-l} = \r^{-l}\Zp-q^{-3/2}[l]_{1/q}z\r^{-l}.
\eq
 From $<\Zp(\zb\r^{-l})> = 0$, $l\geq 1$,
one finds easily the recursion formula
\be
   [l+1]_{q}<\r^{-l}> = [l]_{q}<\r^{-l+1}>,\;\; l\geq1,
\eq
which gives
\be
   <\r^{-l}> = \frac{1}{[l+1]_{q}}<1>,\;\; l\geq 0.
\eq
Similarly
\be
   <\frac{\zb z}{(1+\zb z)^{l}}> = (\frac{1}{[l]_{q}}-\frac{1}{[l+1]_{q}})<1>,
    \;\; l\geq 1.
\eq
We leave it to the reader to find the expression for
\be
   <\frac{(\zb z)^{k}}{(1+\zb z)^{l}}>,\;\; l\geq k.
\eq

The above results can also be obtained by using the relations
(\ref{SU1}) and (\ref{SU2}) for $z$ and $\zb$ in terms of the $SU_{q}(2)$
parameters and known results\cite{W1,Z} for the
Haar measure of $SU_{q}(2)$.
However we wanted to show that one can formulate the integration
directly for the sphere.

As an application of the stereographic projection, we can define an integration
on the complex quantum plane
$C_q$ by inserting an appropriate measure factor
$\rho^2$. $C_q$ has the same algebra (\ref {zz}) and differential calculus
(\ref {zdz}) to (\ref {deldelb}),
but a different $*$-structure (\ref {Cdz*}) to (\ref{Cdelb*}).
Classically, it holds
$\int_C dz d\zb/{2 \pi i} f(z,\zb) =
\int_{S^2} \rho^2 f(z, \zb) $
Motivated by this, we define an integration over the quantum plane as,
\be \label{int-C}
\int_{C_q} f(z,\zb) \equiv <  \rho^2 f(z, \zb) >.
\eq

We need to check that this integration is translationally invariant, namely,
$\int \del f= \int \delb f =0.$
To show this, we must find relations between  the infinitesmal
generators $\del, \delb$ on the plane  and $\Zp, \Zm, \H$ on the sphere.
Introduce the differential operators,
\bee
&C =1-\lam q^{-1} z \del,\\
&D= 1+\lambda q \zb \delb
\eqq
and
\be
B= 1-\lam q^{-1} z \del +\lambda q \zb \delb -\lambda^2 q^{-2} \rho \delb \del.
\eq

One finds the following realizations of $\Zp, \Zm, \H$ as pseudo-differential
operators,
which satisfy Eqs. (\ref{ZZH}) to (\ref{Z*}):
\bee \label{ZZH=pseudo}
&q^{3/2} \Zp = (z^2 \del  +q^2 \delb  B^{-1}) C^{-1},\\
&-q^{3/2} \Zm = (q^2 \zb^2 \delb  + \del B^{-1}) D^{-1}
\eqq
and
\be
\H = \frac{1-B^{-2}}{1-q^2}.
\eq

One also has
\be \label{del=ZZH}
q^{-1} \rho^2 \delb=(\Zm z \Zp -q^4 \Zp z \Zm +q^{1/2}(1+q^2)\Zp) B\\
\eq
and
\be \label{delb=ZZH}
q^{-1} \rho^2 \del=(q^4 \Zp \zb \Zm - \Zm \zb \Zp -q^{1/2}(1+q^2)\Zm) B.
\eq

Together with  the definition (\ref {int-C}), we have
\be
\int_{C_q} \del f = < \rho^2 \del f >=  < \Zp \cdots>-<\Zm \cdots>
\eq
and
\be
\int_{C_q} \delb f= < \rho^2 \del f > = < \Zm \cdots>-<\Zp \cdots> ,
\eq
which are both zero since the integral on the sphere is defined by $< \O f> =0$
for $\O=
{\cal Z_\pm, H}$. So
the integral defined by (\ref{int-C}) is translational invariant.

\section{ACKNOWLEDGEMENTS}

We would like to thank Paul Watts for important discussions in the early
stages of this work and Piotr Podle\'{s} for useful comments.

This work was supported in part by the Director, Office of
Energy Research, Office of High Energy and Nuclear Physics, Division of
High Energy Physics of the U.S. Department of Energy under Contract
DE-AC03-76SF00098 and in part by the National Science Foundation under
grant PHY-90-21139.

\end{document}